\documentstyle[11pt]{article}
\begin{document}
\newcommand{\qed}{\hphantom{.}\hfill $\Box$\medbreak}
\newcommand{\A}{{\cal A}}
\newcommand{\B}{{\cal B}}
\newcommand{\U}{{\cal U}}
\newcommand{\G}{{\cal G}}
\newcommand{\cZ}{{\cal Z}}
\newcommand{\proof}{\noindent{\bf Proof \ }}
\renewcommand{\theequation}{\thesection.\arabic{equation}}
\newtheorem{Theorem}{Theorem}[section]
\newtheorem{Lemma}[Theorem]{Lemma}
\newtheorem{Corollary}[Theorem]{Corollary}
\newtheorem{Remark}[Theorem]{Remark}
\newtheorem{Example}[Theorem]{Example}
\newtheorem{Definition}[Theorem]{Definition}
\newtheorem{Construction}[Theorem]{Construction}
\newcounter{Ictr}
\def\binom#1#2{{#1\choose#2}}
\newcommand{\Item}{\refstepcounter{Ictr}\item[(\theIctr)]}

\thispagestyle{empty}

\title{Decompositions of $3$-uniform hypergraph $K_v^{(3)}$ into
hypergraph $K_4^{(3)}+e$
 \footnote{Supported by the Research Foundation of Beijing Jiaotong University under Grant No.
2008RC036 (T. Feng), and by NSFC grant 10771013 (Y. Chang).}}
\author{Tao Feng,  Yanxun Chang\\
Institute of Mathematics \\
Beijing Jiaotong University\\
Beijing 100044, P. R. China\\
tfeng@bjtu.edu.cn\\
yxchang@bjtu.edu.cn\\}
\date{}
\maketitle

\noindent {\bf Abstract:} In this paper it is established that a
decomposition of a $3$-uniform hypergraph $K_v^{(3)}$ into a special
kind of hypergraph $K_4^{(3)}+e$ exists if and only if $v\equiv
0,1,2$ $({\rm mod}$ $5)$ and $v\geq 7$.



\section{Introduction}

A {\em hypergraph} $H$ is a pair $(V,E)$, where $V$ is a finite set
of vertices, $E$ is a family of subsets of $V$ (called {\em
hyperedges} or {\em edges}). A hypergraph is called {\em simple} if
$E$ has no repeated edges. All hypergraphs considered in this paper
are simple. A {\em sub-hypergraph} $H'=(V',E')$ of $H=(V,E)$ is a
hypergraph satisfying $V'\subseteq V$ and $E'\subseteq E$.

A hypergraph is said to be {\em $t$-uniform} if each of its edges
contains exactly $t$ vertices. In particular a $2$-uniform
hypergraph is just a graph. For a $t$-uniform hypergraph $H$, let
$V(H)$ and $E(H)$ denote the vertex-set and edge-set of $H$,
respectively. We say that $H$ contains a vertex $x$ if $x\in V(H)$,
and $H$ contains a set of vertices $\{x_1,x_2,\ldots,x_t\}$ if
$\{x_1,x_2,\ldots,x_t\}\in E(H)$. A $t$-uniform hypergraph is said
to be {\em complete} if the edge-set $E$ contains each $t$-subset of
$V$ exactly once. It is denoted by $K_v^{(t)}$, where $v=|V|$ is
called the {\em order} of the $t$-uniform hypergraph. The {\em
degree} of a vertex $x$ in a hypergraph is the number of edges that
contain this vertex. It is denoted by $d(x)$. For more information
on hypergraphs, the reader may refer to \cite{berge}.

Let $H$ be a $t$-uniform hypergraph and $\Gamma$ be a set of
$t$-uniform hypergraphs. A {\em decomposition} of $H$ into
hypergraphs of $\Gamma$ is a partition of the edges of $H$ into
sub-hypergraphs each of which is isomorphic to a hypergraph in
$\Gamma$. Such a decomposition of $H$ into $\Gamma$ is denoted by
$(H,\Gamma)$-design. Hypergraph decompositions have an interesting
application in secret sharing schemes (cf. \cite{cgal}). When
$H=K_v^{(t)}$, a $(K_v^{(t)},\Gamma)$-design is called a {\em
$t$-wise balanced ${\it \Gamma}$ design}, denoted by
$S(t,\Gamma,v)$.  If $\Gamma$ only contains one hypergraph $J$, we
write $S(t,\{J\},v)$ simply as $S(t,J,v)$. Let $K$ be a set of
positive integers and $\Omega$ a set of complete $t$-uniform
hypergraphs, where the order of each element in $\Omega$ is from
$K$. We denote an $S(t,\Omega,v)$ by $S(t,K,v)$, which corresponds
to the traditional concept of $t$-wise balanced design ($t$-BD)
\cite{bjl}. Therefore the $t$-wise balanced $\Gamma$ design is a
generalization of the $t$-wise balanced design.

One of the interesting problems in design theory is to determine the
existence spectrum of $S(t,\Gamma,v)$, where $\Gamma$ is a set of
$t$-uniform hypergraphs. When $t=2$, significant progress was made
on this problem by many authors (e.g., see \cite{be} and the
references therein). However much less is known about $t\geq 3$. For
$t=3$, the necessary conditions for the existence of an
$S(3,\Gamma,v)$ are as follows.

\begin{Lemma}
\label{necessary} $({\rm \cite{fc}})$ Let $\Gamma$ be any set of
$3$-uniform hypergraphs and $J\in \Gamma$. Let $d_J^{*}(x_1,x_2)$ be
the number of edges in $J$ containing the two vertices $x_1$ and
$x_2$. The following are necessary conditions for the existence of
an $S(3,\Gamma,v)$:

\begin{enumerate}
\item[{\rm (1)}] $v\geq min\{|V(J)|:J \in \Gamma\}$;

\item[{\rm (2)}] $\binom v3 \equiv 0$ $({\rm mod }$ $d_0)$, where
$d_0=gcd\{|E(J)|:J \in \Gamma\}$;

\item[{\rm (3)}] $\binom {v-1} 2 \equiv 0$ $({\rm mod }$ $d_1)$,
where $d_1=gcd\{d(x):x\in V(J),J\in \Gamma\}$;

\item[{\rm (4)}] $\binom {v-2} 1 \equiv 0$ $({\rm mod }$ $d_2)$,
where $d_2=gcd\{d_J^{*}(x_1,x_2):x_1\neq x_2$, $x_1$, $x_2$ $\in
V(J),J\in\Gamma\}$.
\end{enumerate}
\end{Lemma}

For an edge $e\in E(K_4^{(3)})$, let $K_4^{(3)}-e$ denote the
hypergraph obtained from $K_4^{(3)}$ by deleting the edge $e$. In
\cite{fc}, the authors investigated the existence of an
$S(3,K_4^{(3)}-e,v)$ as follows.

\begin{Theorem}
\label{K_4-e} $({\rm \cite{fc}})$ An $S(3,K_4^{(3)}-e,v)$ exists if
and only if $v\equiv 0,1,2$ $({\rm mod}$ $9)$ and $v\geq 9$.
\end{Theorem}

In \cite{hanani60} Hanani gave the following theorem.

\begin{Theorem}
\label{K_4} $({\rm \cite{hanani60}})$ An $S(3,K_4^{(3)},v)$ exists
if and only if $v\equiv 2,4$ $({\rm mod}$ $6)$ and $v\geq 4$.
\end{Theorem}

Combining Theorems \ref{K_4-e} and \ref{K_4}, it is natural to
consider what the existence spectrum for an $S(3,K_4^{(3)}+e,v)$ is,
where $K_4^{(3)}+e$ is a hypergraph $(V,E)$ with $V=\{1,2,3,4,5\}$
and $E=\{\{1,2,3\}$, $\{1,2,4\},$ $\{1,3,4\}$, $\{2,3,4\}$,
$\{3,4,5\}\}$. In this paper we establish that

\begin{Theorem}
\label{K_4+e}  An $S(3,K_4^{(3)}+e,v)$ exists if and only if
$v\equiv 0,1,2$ $({\rm mod}$ $5)$ and $v\geq 7$.
\end{Theorem}

For convenience, in this paper we always assume that $K$ is a set of
positive integers, $\Gamma$ is a set of $t$-uniform hypergraphs, and
$\Omega$ is a set of complete $t$-uniform hypergraphs, where the
order of each element in $\Omega$ is from $K$.

\section{Recursive constructions}

To describe our recursive constructions, we need the following
auxiliary designs. For the general background on design theory, the
reader is referred to \cite{bjl}.

Let $n$ and $t$ be positive integers. Suppose that $X$ is a set of
points, ${\cal B}$ is a collection of hypergraphs on the subsets of
$X$ (called {\em blocks}), and ${\cal G}$ is a partition of $X$ into
$n$ non-empty subsets (called {\em groups} or {\em holes}). A {\em
group divisible $({\it \Gamma},t)$-design} is a triple $(X,{\cal
G},{\cal B})$, where for each $B\in {\cal B}$, $B$ is isomorphic to
a hypergraph in $\Gamma$, such that each edge from the edge-set of
each block intersects any given group in at most one point, and each
$t$-subset of $X$ from $t$ distinct groups is contained in a unique
block.

We use the usual exponential notation for the types of group
divisible $(\Gamma,t)$-designs. Then type $g_1^{a_1} g_2^{a_2}
\cdots g_m^{a_m}$ denotes that there are $a_i$ groups of size $g_i$,
$1 \leq i \leq m$. For brevity, a group divisible
$(\Gamma,t)$-design of type $g_1^{a_1} g_2^{a_2} \cdots g_m^{a_m}$
can be denoted by GDD$(t,\Gamma,v)$ of type $g_1^{a_1} g_2^{a_2}
\cdots g_m^{a_m}$, where $v=\sum_{i=1}^{m}a_i g_i$. If $\Gamma$
contains only one hypergraph $J$, we write GDD$(t,\{J\},v)$ as
GDD$(t,J,v)$.

If we replace $\Gamma$ by $\Omega$, then a GDD$(t,\Omega,v)$ is
denoted by GDD$(t,K,$ $v)$, which corresponds to the traditional
concept of group divisible $t$-design ($t$-GDD) \cite{mr}. Therefore
the group divisible $(\Gamma,t)$-design is a generalization of the
group divisible $t$-design. Furthermore, if all the $n$ groups have
the same size $g$, a GDD$(t,K,v)$ is called an $H$ design
\cite{m74}, denoted by $H(n,g,K,t)$.

\begin{Lemma}
\label{GDD-46}$({\rm \cite{ji}})$ For any integer $4\leq n\leq 27$
and $n\neq 5,21$, there is a GDD$(3,\{4,6\},2n)$ of type $2^n$.
\end{Lemma}

The following construction is a variation of the fundamental
construction for $t$-GDD \cite{mr}.

\begin{Construction}
\label{GDD-FC} $({\rm \cite{fc}})$ Suppose that there exists a
GDD$(t,K,v)$ of type $g_1^{a_1} g_2^{a_2} \cdots g_m^{a_m}$. If
there exists a GDD$(t,\Gamma,hk)$ of type $h^k$ for each $k\in K$,
then there exists a GDD$(t,\Gamma,hv)$ of type $(hg_1)^{a_1}
(hg_2)^{a_2} \cdots $ $(hg_m)^{a_m}$.
\end{Construction}

Let $v$, $m$  and $t$ be positive integers, and $s$ be a
non-negative integer. Suppose that $X$ is a set of $v=s+
\sum_{_{1\leq i \leq m}} a_i g_i$ points, $S$ is a subset of $X$ of
size $s$ (called {\em stem}), ${\cal T}=\{G_1,G_2,\ldots,G_n\}$ is a
partition of $X \setminus S$ of type $g_1^{a_1} \cdots g_m^{a_m}$,
$n=\sum_{i=1}^m a_i$, (called {\em groups} or {\em branches}), and
${\cal A}$ is a collection of hypergraphs on the subsets of $X$
(called {\em blocks}). A {\em candelabra $({\it \Gamma},t)$-system}
is a quadruple $(X,S,{\cal T},{\cal A})$ of type $(g_1^{a_1} \cdots
g_m^{a_m}:s)$, where for each $A\in {\cal A}$, $A$ is isomorphic to
one of $\Gamma$, such that every $t$-subset $T\subset X$ with
$|T\cap (S\cup G_i)|<t$ for all $i$ is contained in a unique block
and no $t$-subset of $S\cup G_i$ is contained in any block. Such a
system is denoted by CS$(t,\Gamma,v)$ of type $(g_1^{a_1} \cdots
g_m^{a_m}:s)$.  If $\Gamma$ contains only one hypergraph $J$, we
write CS$(t,\{J\},v)$ as CS$(t,J,v)$.

If we replace $\Gamma$ by $\Omega$, then a CS$(t,\Gamma,v)$ is
denoted by CS$(t,K,v)$, which corresponds to the traditional concept
of candelabra $t$-design \cite{mr}. Thus the candelabra $(
\Gamma,t)$-system is a generalization of the candelabra $t$-system.

\begin{Lemma}
\label{2^n:2} There exists a $CS(3,\{4,6\},2n+2)$ of type $(2^n:2)$
for any integer $n\geq 3$.
\end{Lemma}

\proof For $n\equiv 0,1$ $({\rm mod}$ $3)$ and $n\geq 3$, let
$(X,{\cal B})$ be an $S(3,4,2(n+1))$ by Theorem \ref{K_4}. Let $a$
and $b$ be two distinct points in $X$. Let ${\cal T}=\{B\setminus
\{a,b\}:B\in {\cal B},\{a,b\}\subset B\}$ and ${\cal D}=\{B:B\in
{\cal B},\{a,b\}\subset B\}$. Then it is readily checked that
$(X,\{a,b\},{\cal T},{\cal B}\setminus {\cal D})$ is a
$CS(3,4,2n+2)$ of type $(2^n:2)$.

For $n\equiv 2$ $({\rm mod}$ $3)$ and $n\geq 5$, let
$(X',\emptyset,{\cal G},{\cal B'})$ be a $CS(3,4,2(n+1))$ of type
$(6^{(n+1)/3}:0)$ by Theorem $1$ in \cite{m74}. Let $c$ and $d$ be
two distinct points from two distinct groups in $\cal G$. Let ${\cal
T'}=\{B\setminus \{c,d\}:B\in {\cal B'},\{c,d\}\subset B\}$ and
${\cal D'}=\{B:B\in {\cal B'},\{c,d\}\subset B\}$. Then it is
readily checked that $(X',\{c,d\},{\cal T'},({\cal B'}\setminus
{\cal D'})\cup {\cal G})$ is a $CS(3,\{4,6\},2n+2)$ of type
$(2^n:2)$. This completes the proof. \qed

By Lemma \ref{2^n:2}, the following lemma is straightforward.

\begin{Lemma}
\label{46} $({\rm \cite{hanani63}})$ There exists an
$S(3,\{4,6\},v)$ for any integer $v\equiv 0$ $({\rm mod}$ $2)$ and
$v\geq 4$.
\end{Lemma}

We quote the following result for later use.

\begin{Lemma}
\label{45679} $({\rm \cite{hanani63,ji45}})$ There exists an
$S(3,\{4,5,6,7,9,11,13,15,19,$ $23,27\},v)$ for any integer $v\geq
4$.
\end{Lemma}

\begin{Construction}
\label{g^2-construction} $({\rm \cite{fc}})$ Suppose that there
exists a $GDD(3,\Gamma,gn)$ of type $g^n$. If there is a
$CS(3,\Gamma,2g+s)$ of type $(g^2:s)$, then there is a
$CS(3,\Gamma,gn+s)$ of type $(g^n:s)$.
\end{Construction}

Let $s\leq v$ be a non-negative integer. For convenience, in what
follows a $(K_v^{(t)}\setminus K_s^{(t)},\Gamma)$-design is denoted
by HS$(t,\Gamma;v,s)$, where the set of the $s$ points is called the
{\em hole} of this design. The following construction is simple but
useful.

\begin{Construction}
\label{filling construction} Suppose that there exists a
$CS(t,\Gamma,v)$ of type $(g_1^{a_1} \cdots
g_{m-1}^{a_{m-1}}g_m^1:s)$. If there is an $HS(t,\Gamma;g_i+s,s)$
for each $1\leq i\leq m-1$ and there is an $S(t,\Gamma,g_m+s)$, then
there is an $S(t,\Gamma,v)$.
\end{Construction}

\begin{Construction}
\label{HFC} Let $(X,S,{\cal G},{\cal A})$ be a $CS(3,K,v)$ of type
$(g_1^{a_1} \cdots $ $g_m^{a_m}:s)$ with $S=\{x_1,x_2,\ldots,x_s\}$.
Suppose that for each block $A\in {\cal A}$ containing $x_1$, there
exists a $CS(3,\Gamma,$ $b(|V(A)|-1)+r)$ of type $(b^{|V(A)|-1}:r)$.
Suppose that for each block $A\in {\cal A}$ containing $x_i$ for
$2\leq i\leq s$, there exists a $GDD(3,\Gamma,b|V(A)|)$ of type
$b^{|V(A)|}$. If for each block $A\in {\cal A}$ not containing $x_i$
for any $1\leq i\leq s$, there exists a $GDD(3,\Gamma,b|V(A)|)$ of
type $b^{|V(A)|}$, then there exists a $CS(3,\Gamma,v')$ of type
$((bg_1)^{a_1} \cdots (bg_m)^{a_m}:r+sb-b)$, where $v'=(v-1)b+r$.
\end{Construction}

\proof For each block $A\in {\cal A}$ containing $x_1$, construct a
$CS(3,\Gamma,$ $b(|V(A)|-1)+r)$ of type $(b^{|V(A)|-1}:r)$ on
$((V(A)\setminus \{x_1\})\times Z_b) \bigcup S'$ with groups
$\{x\}\times Z_b$, $x\in V(A) \setminus \{x_1\}$ and a stem $S'$ of
size $r$. We denote its block set by ${\cal B}_{A(x_1)}$. For each
block $A\in {\cal A}$ containing $x_i$, $2\leq i\leq s$, construct a
$GDD(3,\Gamma,b|V(A)|)$ of type $b^{|V(A)|}$ on $V(A)\times Z_b$
with groups $\{x\}\times Z_b$, $x\in V(A)$. We denote its block set
by ${\cal B}_{A(x_i)}$. For each block $A\in {\cal A}$ not
containing $x_i$ for any $1\leq i\leq s$, construct a
$GDD(3,\Gamma,b|V(A)|)$ of type $b^{|V(A)|}$ on $V(A)\times Z_b$
with groups $\{\{x\}\times Z_b: x\in V(A)\}$. We denote its block
set by ${\cal B}_A$.

Let $X'=((X\setminus \{x_1\}) \times Z_b) \bigcup S'$ and ${\cal
G}'=\{G \times Z_b : G\in {\cal G}\}$. Let $S''=
(\cup_{i=2}^s(\{x_i\}\times Z_b))\bigcup S'$. For $1\leq i\leq s$,
let ${\cal B}_i$ denote the union of ${\cal B}_{A(x_i)}$ for all
blocks $A\in {\cal A}$ containing $x_i$. Let ${\cal B}'$ denote the
union of ${\cal B}_A$ for all blocks $A\in {\cal A}$ not containing
$x_i$ for any $1\leq i\leq s$. Let ${\cal B}=(\cup_{i=1}^s {\cal
B}_i) \cup {\cal B}'$. Then using similar arguments as in
Construction $2.8$ in \cite{fc}, it is readily checked that
$(X',S'',{\cal G}',{\cal B})$ is the required $CS(3,\Gamma,v')$.
\qed

Construction $2.8$ in \cite{fc} can be seen as a corollary of the
above construction.

\section{Direct constructions}

In the following we always denote the copy of $K_4^{(3)}+e$ with
vertices $x$, $y$, $z$, $u$, $v$ and edges $\{x,y,z\}$, $\{x,y,u\}$,
$\{x,z,u\}$, $\{y,z,u\}$, $\{z,u,v\}$  by $(x,y,z,u,v)$.

\begin{Lemma}
\label{71116} There exists an $S(3,K_4^{(3)}+e,v)$ for $v \in
\{7,11,16,26,31$, $32\}$.
\end{Lemma}

\proof Let $X=Z_v$. Base blocks for these designs are given below.
All other blocks are obtained by developing these base blocks by
$+1$ modulo $v$.

\begin{tabular}{llll}
$v=7$: &$(0,1,2,4,5)$.\\
$v=11$:& $(0,1,2,6,4)$,& $(0,1,3,9,7)$,& $(0,1,4,7,8)$.\\
$v=16$:& $(0,1,2,4,5)$,& $(0,1,5,6,12)$,& $(0,1,8,10,3)$,\\
& $(0,1,9,13,15)$, & $(0,2,5,10,13)$,& $(0,2,6,13,8)$,\\
& $(0,3,6,10,15)$. \\
$v=26$:& $(0,1,2,4,5)$,& $(0,1,5,6,2)$,& $(0,1,7,8,16)$,\\
& $(0,1,10,11,22)$,& $(0,1,13,15,2)$,& $(0,1,14,18,2)$,\\&
$(0,2,5,7,11)$,& $(0,2,8,10,4)$,& $(0,2,9,12,3)$,\\&
$(0,2,11,17,3)$,& $(0,2,16,19,3)$,& $(0,3,6,21,10)$,\\&
$(0,3,7,11,24)$,& $(0,3,8,15,22)$,& $(0,3,9,14,2)$,\\&
$(0,3,16,22,4)$,& $(0,4,9,18,2)$,& $(0,4,10,16,3)$,\\&
$(0,4,11,21,6)$,& $(0,5,10,18,3)$.\\
$v=31$:& $(0,1,2,4,5)$,& $(0,1,5,6,2)$,& $(0,1,7,8,16)$,\\&
$(0,1,10,11,2)$, & $(0,1,12,13,27)$,& $(0,1,14,17,2)$,\\&
$(0,1,16,18,2)$,& $(0,2,5,7,11)$,& $(0,2,8,10,4)$,\\&
$(0,2,9,11,21)$,& $(0,2,13,19,3)$,& $(0,2,14,20,5)$,\\&
$(0,2,18,21,3)$,& $(0,3,6,10,2)$,& $(0,3,8,11,2)$,\\&
$(0,3,9,14,2)$,& $(0,3,12,21,5)$,& $(0,3,20,24,5)$,
\end{tabular}

\begin{tabular}{llll}
& $(0,3,22,27,17)$,& $(0,4,9,21,15)$,& $(0,4,10,20,3)$,\\&
$(0,4,12,19,5)$,& $(0,4,13,17,5)$,& $(0,4,15,24,10)$,\\&
$(0,5,11,19,28)$,& $(0,5,13,23,30)$,& $(0,5,15,20,2)$,\\&
$(0,6,13,24,5)$,& $(0,6,15,23,8)$.\\
$v=32$:& $(0,1,2,4,5)$,& $(0,1,5,6,2)$,& $(0,1,7,8,16)$,\\&
$(0,1,10,11,2)$,& $(0,1,12,13,26)$,& $(0,1,15,16,2)$,\\&
$(0,2,5,7,11)$,& $(0,2,8,10,4)$,& $(0,2,9,11,21)$,\\&
$(0,2,13,15,3)$,& $(0,2,14,16,31)$,& $(0,3,6,10,2)$,\\&
$(0,3,8,11,2)$,& $(0,3,9,14,2)$,& $(0,3,12,15,2)$,\\&
$(0,3,13,17,2)$,& $(0,3,16,21,2)$,& $(0,3,18,25,2)$,\\&
$(0,3,19,28,2)$,& $(0,4,9,18,5)$,& $(0,4,10,16,29)$,\\&
$(0,4,11,25,15)$,& $(0,4,12,20,2)$,& $(0,4,15,27,5)$,\\&
$(0,4,17,24,5)$,& $(0,5,10,17,4)$,& $(0,5,11,21,3)$,\\
&$(0,5,13,22,2)$,& $(0,5,15,24,9)$,& $(0,6,14,20,3)$,\\&
$(0,6,17,25,4)$.
\end{tabular} \\\qed

\begin{Lemma}
\label{101215} There exists an $S(3,K_4^{(3)}+e,v)$ for $v \in
\{10,12,15\}$.
\end{Lemma}
\proof  For $v=10$: let $X=Z_{10}$. All $24$ blocks are listed
below.

\hspace*{-0.4 cm}
\begin{tabular}{llll}
$(0,1,8,9,5)$,& $(0,1,4,5,6)$,& $(0,1,6,7,4)$,& $(0,2,5,8,4)$,\\
$(0,2,7,9,5)$,& $(0,4,3,9,8)$,& $(0,5,3,7,4)$,& $(0,3,6,8,5)$,\\
$(0,8,4,7,9)$,& $(0,6,5,9,4)$,& $(4,7,1,2,0)$,& $(1,2,5,9,3)$,\\
$(1,2,6,8,7)$,& $(4,8,1,3,0)$,& $(1,6,3,5,8)$,& $(7,9,1,3,2)$,\\
$(1,9,4,6,8)$,& $(8,1,5,7,4)$,& $(4,5,2,3,0)$,& $(2,9,3,6,7)$,\\
$(2,3,7,8,9)$,& $(2,4,8,9,6)$,& $(0,2,4,6,3)$,& $(2,5,6,7,9)$.
\end{tabular}

For $v=12$: let $X=Z_{11}\cup \{\infty\}$. Base blocks for this
design are given below. All other blocks are obtained by developing
these base blocks by $+1$ modulo $11$, where $\infty+1=\infty$.

\begin{tabular}{llll}
$(8,1,0,5,\infty)$,& $(0,1,\infty,3,7)$,& $(1,4,6,0,2)$,&
$(7,9,1,0,2)$.
\end{tabular}

For $v=15$: let $X=Z_{13}\cup \{\infty_1,\infty_2\}$. Base blocks
for this design are given below. All other blocks are obtained by
developing these base blocks by $+1$ modulo $13$, where
$\infty_i+1=\infty_i$ for $i=1,2$.

\begin{tabular}{lll}
$(0,1,4,\infty_1,\infty_2)$,& $(\infty_1,7,0,2,4)$,&
$(\infty_2,5,0,7,1)$,\\ $(\infty_2,10,0,1,9)$,& $(2,6,0,1,3)$,&
$(1,11,0,8,3)$,\\
\end{tabular}

\begin{tabular}{lll}
 $(2,5,0,9,3)$.
\end{tabular}\qed

\begin{Lemma}
\label{GDD-5^4-5^6-10^5} There exists a $GDD(3,K_4^{(3)}+e,gn)$ of
type $g^n$ for $(g,n)\in \{(5,4),(5,6),(10,5)\}$.
\end{Lemma}
\proof For $(g,n)=(5,4)$: let $X=Z_{20}$ and ${\cal G}=\{4Z_5+j:
0\leq j\leq 3\}$. Base blocks for this design are given below. All
other blocks are obtained by developing these base blocks by $+1$
modulo $20$.

\hspace*{-0.3 cm}
\begin{tabular}{llll}
$(0,1,2,7,12)$,& $(0,1,3,14,4)$,&
$(0,1,15,18,5)$,& $(0,3,9,18,8)$,\\ $(0,5,7,14,17)$. \\
\end{tabular}

For $(g,n)=(5,6)$: let $X=Z_{25} \cup
\{\infty_0,\infty_1,\ldots,\infty_4\}$ and ${\cal G}=\{5Z_5+i: 0\leq
i\leq 4\}\cup \{\infty_0,\infty_1,\ldots,\infty_4\}$. Base blocks
for this design are given below. All other blocks are obtained by
developing these base blocks by $+1$ modulo $25$, where
$\infty_j+1=\infty_{j+1}$, $0\leq j\leq 4$, the subscripts are
reduced modulo $5$.

\begin{tabular}{lll} $(\infty_0,2,0,1,9)$,& $(\infty_0,3,0,4,11)$,&
$(\infty_0,8,0,6,18)$,\\
$(\infty_0,9,0,7,14)$,& $(\infty_0,14,0,11,8)$,&
$(\infty_0,12,0,13,19)$,\\ $(0,16,22,\infty_0,1)$,&
$(0,17,24,\infty_0,1)$,& $(0,18,21,\infty_0,3)$,\\
$(0,19,23,\infty_0,2)$,& $(4,11,23,\infty_0,7)$,&
$(3,14,22,\infty_0,8)$,\\
 $(1,14,17,\infty_0,4)$,&
$(4,13,21,\infty_0,7)$,& $(0,1,4,7,5)$,\\
 $(0,1,8,12,4)$,&
$(0,1,14,23,2)$,& $(0,2,4,13,7)$,\\
\end{tabular}

\begin{tabular}{lll}
 $(0,2,6,14,3)$,&
$(0,6,7,23,14)$.
\end{tabular}

For $(g,n)=(10,5)$: let $X=Z_{50}$ and ${\cal G}=\{5Z_{10}+j: 0\leq
j\leq 4\}$. The $40$ base blocks for this design can be obtained by
multiplying each of the following $8$ base blocks by $(11)^i$,
$i=0,1,2,3,4$. All other blocks are obtained by developing these
base blocks by $+1$ modulo $50$.

\begin{tabular}{lll}
$(0,1,2,4,5)$,& $(0,1,7,8,4)$,& $(0,1,9,33,5)$,\\ $(0,1,13,19,5)$,&
$(0,1,17,23,10)$,& $(0,1,18,42,5)$,\\ $(0,1,27,29,3)$,&
$(0,1,28,34,41)$.
\end{tabular}
\\\qed

\begin{Lemma}
\label{g^n:0} There exists a $CS(3,K_4^{(3)}+e,gn)$ of type
$(g^n:0)$ for $(g,n)\in \{(5,3)$, $(5,5)$, $(10,2)$, $(15,2)\}$.
\end{Lemma}

\proof Let $X=Z_{gn}$, $S=\emptyset$ and ${\cal T}=\{nZ_g+j:0\leq
j\leq n-1\}$. Base blocks for these designs are given below.

For $(g,n)=(5,3)$, develop the following base blocks by $+3$ modulo
$15$.

\tabcolsep 0.03in
\begin{tabular}{llll}
$(0,6,1,7,14)$,& $(1,8,0,11,4)$,& $(0,4,2,6,11)$,& $(0,5,2,7,10)$,\\
$(0,8,2,9,6)$,& $(0,2,10,11,5)$,& $(2,5,1,4,3)$,& $(1,2,0,3,5)$,\\
$(1,4,0,12,11)$,& $(1,5,0,14,2)$,& $(3,13,0,10,14)$,&
$(4,13,0,5,12)$,\\
$(4,10,0,8,5)$,& $(0,4,9,14,3)$,& $(2,7,1,9,3)$,& $(2,3,10,14,4)$,\\
$(1,14,4,8,6)$.
\end{tabular}

For $(g,n)=(5,5)$, develop the following base blocks by $+1$ modulo
$25$.

\tabcolsep 0.03in
\begin{tabular}{llll}
$(0,1,2,4,5)$,& $(0,1,5,6,2)$,& $(0,1,7,8,16)$,& $(0,1,10,12,2)$,\\
$(0,1,11,14,2)$,& $(0,1,13,16,5)$,& $(0,1,15,17,2)$,&
$(0,2,5,7,11)$,\\ $(0,2,8,13,20)$,& $(0,2,9,19,16)$,&
$(0,2,14,21,23)$,& $(0,3,6,10,13)$,\\ $(0,3,8,19,11)$,&
$(0,3,9,20,24)$,& $(0,4,8,16,3)$,& $(0,4,9,15,3)$,\\
$(0,4,10,17,6)$,& $(0,4,13,20,2)$.
\end{tabular}

For $(g,n)=(10,2)$, develop the following base blocks by $+1$ modulo
$20$.

\tabcolsep 0.03in
\begin{tabular}{llll}
$(0,1,2,5,3)$,& $(0,1,6,7,14)$,& $(0,1,9,16,2)$,& $(0,1,10,13,7)$,\\
$(0,1,11,18,14)$,& $(0,1,12,17,10)$,& $(0,2,9,17,12)$,&
$(0,2,11,15,18)$,\\ $(0,3,9,14,4)$.
\end{tabular}

For $(g,n)=(15,2)$, develop the following base blocks by $+1$ modulo
$30$.

\tabcolsep 0.01in
\begin{tabular}{llll}
$(0,1,2,5,3)$,& $(0,1,6,7,2)$,& $(0,1,8,9,5)$,& $(0,1,10,11,8)$,\\
$(0,1,12,15,2)$,& $(0,1,13,18,2)$,& $(0,1,14,17,5)$,&
$(0,1,16,19,8)$,\\ $(0,2,7,9,4)$,& $(0,2,11,19,16)$,&
$(0,2,13,17,4)$,& $(0,2,15,21,10)$,\\
\end{tabular}

\tabcolsep 0.01in
\begin{tabular}{llll}
 $(0,3,6,11,16)$,&
$(0,3,7,12,19)$,& $(0,3,9,26,4)$,& $(0,3,10,23,15)$,\\
$(0,3,13,24,8)$,& $(0,4,11,23,2)$,& $(0,4,15,25,10)$,&
$(0,5,14,21,7)$,\\ $(0,6,13,21,28)$.
\end{tabular}
\\\qed

\begin{Lemma}
\label{g^n:1} There exists a $CS(3,K_4^{(3)}+e,gn+1)$ of type
$(g^n:1)$ for $(g,n)\in \{(5,3)$, $(5,5)$, $(6,4)$, $(10,2)\}$.
\end{Lemma}

\proof Let $X=Z_{gn}\cup \{\infty\}$, $S=\{\infty\}$ and ${\cal
T}=\{nZ_g+j:0\leq j\leq n-1\}$. Base blocks for these designs are
given below.

For $(g,n)=(5,3)$, develop the following base blocks by $+3$ modulo
$15$, where $\infty+3=\infty$.

\begin{tabular}{llll}
$(\infty,2,1,0,9)$,& $(0,4,8,\infty,9)$,& $(0,5,10,\infty,2)$,&
$(0,11,13,\infty,6)$,\\ $(0,7,1,5,10)$,& $(0,11,1,8,10)$,&
$(0,4,2,5,13)$,& $(0,1,3,4,2)$,\\ $(0,6,1,10,2)$,& $(0,7,2,3,10)$,&
$(2,8,0,6,5)$,& $(0,12,2,10,11)$,\\ $(0,14,2,13,3)$,&
$(0,3,8,14,4)$,& $(0,11,3,10,13)$,& $(0,6,4,13,14)$,\\
$(0,11,4,7,8)$,& $(0,5,8,9,14)$,& $(1,3,11,13,8)$,& $(1,11,2,4,8)$.
\end{tabular}

For $(g,n)=(5,5)$, develop the following base blocks by $+1$ modulo
$25$, where $\infty+1=\infty$.

\begin{tabular}{llll}
$(1,8,0,9,\infty)$,& $(10,11,0,1,\infty)$,& $(1,14,0,12,\infty)$,&
$(1,13,0,19,\infty)$,\\ $(1,20,0,21,\infty)$,&
$(1,24,0,22,\infty)$,& $(2,5,0,7,\infty)$,& $(2,6,0,8,\infty)$,\\
$(2,9,0,11,\infty)$,& $(10,17,0,2,\infty)$,& $(0,2,12,15,13)$,&
$(0,3,6,10,7)$,\\
$(0,3,8,11,2)$,& $(0,3,9,18,13)$,& $(0,3,12,16,2)$,&
$(0,3,14,21,9)$,\\ $(0,4,8,14,7)$,& $(0,4,11,17,2)$,&
$(0,4,12,20,7)$,& $(0,5,11,16,4)$.
\end{tabular}

For $(g,n)=(6,4)$, develop the following base blocks by $+1$ modulo
$24$, where $\infty+1=\infty$.

\begin{tabular}{llll}
$(2,4,0,1,\infty)$,& $(1,6,0,5,\infty)$,& $(0,8,1,7,\infty)$,&
$(1,10,0,9,\infty)$,\\ $(0,1,11,13,\infty)$,& $(1,12,0,21,\infty)$,&
$(0,22,1,14,\infty)$,& $(2,5,0,7,\infty)$,\\ $(0,2,6,20,\infty)$,&
$(0,2,8,18,6)$,& $(0,2,9,11,19)$,& $(0,3,6,12,21)$,\\
$(0,3,7,20,11)$,& $(0,3,8,19,11)$,& $(0,3,10,13,4)$, &
$(0,4,9,19,5)$,\\ $(0,5,10,17,2)$,& $(0,5,11,18,2)$.
\end{tabular}

For $(g,n)=(10,2)$, develop the following base blocks by $+1$ modulo
$20$, where $\infty+1=\infty$.

\hspace*{-0.4 cm}
\begin{tabular}{llll}
$(2,9,0,1,\infty)$,& $(1,12,0,3,\infty)$,& $(1,4,0,7,\infty)$,&
$(1,16,0,5,\infty)$,\\ $(0,1,6,15,\infty)$,& $(0,1,10,17,12)$,&
$(0,1,11,14,16)$,& $(0,1,13,18,3)$,\\ $(0,2,5,13,8)$,&
$(0,3,7,14,2)$.
\end{tabular}
\\\qed

\begin{Lemma}
\label{g^n:2} There exists a $CS(3,K_4^{(3)}+e,gn+2)$ of type
$(g^n:2)$ for each $(g,n)\in \{(5,3)$, $(5,5)$, $(10,2)\}$.
\end{Lemma}

\proof Let $X=Z_{gn}\cup \{\infty_1,\infty_2\}$,
$S=\{\infty_1,\infty_2\}$ and ${\cal T}=\{nZ_g+j:0\leq j\leq n-1\}$.
Base blocks for these designs are given below.

For $(g,n)=(5,3)$, develop the following base blocks by $+3$ modulo
$15$, where $\infty_i+3=\infty_i$ for $i=1,2$.

\hspace*{-0.5 cm} \tabcolsep 0.03in
\begin{tabular}{llll}
$(\infty_1,1,0,2,13)$,& $(0,4,8,\infty_1,9)$,&
$(0,5,10,\infty_1,2)$,& $(0,11,13,\infty_1,6)$,\\
$(\infty_2,1,0,5,13)$,& $(2,9,7,\infty_2,8)$,&
$(0,4,2,\infty_2,6)$,& $(2,10,3,\infty_2,13)$,\\ $(1,8,0,10,13)$,&
$(1,4,0,3,11)$,& $(0,7,1,9,6)$,& $(0,1,11,14,6)$,\\
$(0,5,2,3,11)$,& $(0,8,2,6,14)$,& $(0,2,7,10,5)$,& $(0,7,3,8,10)$,\\
$(0,14,3,13,12)$,& $(6,10,0,4,14)$,& $(0,4,7,11,14)$,&
$(0,11,5,6,12)$,\\ $(1,5,2,4,11)$,& $(1,8,2,7,14)$,&
$(1,5,7,11,12)$.
\end{tabular}

For $(g,n)=(5,5)$, develop the following base blocks by $+1$ modulo
$25$, where $\infty_i+1=\infty_i$ for $i=1,2$.

\hspace*{-0.8 cm} \tabcolsep 0.01in
\begin{tabular}{llll}
$(0,9,1,8,\infty_1)$,& $(\infty_1,1,3,0,\infty_2)$,&
$(\infty_1,0,13,4,\infty_2)$,& $(\infty_1,6,14,0,\infty_2)$,\\
$(0,4,17,\infty_2,19)$,& $(\infty_2,7,1,0,6)$,& $(0,1,10,11,4)$,&
$(0,1,12,14,2)$,\\ $(0,1,13,20,2)$,& $(0,1,21,24,23)$,&
$(0,2,4,7,3)$,& $(0,2,6,8,15)$,\\ $(0,2,10,17,15)$,&
$(0,2,11,16,5)$,& $(0,2,12,18,3)$,& $(0,3,6,19,23)$,\\
$(0,3,7,11,2)$,& $(0,3,8,14,5)$,& $(0,3,10,13,2)$,&
$(0,3,12,20,7)$,\\ $(0,4,9,16,5)$,& $(0,4,10,14,9)$.
\end{tabular}

For $(g,n)=(10,2)$, develop the following base blocks by $+1$ modulo
$20$, where $\infty_i+1=\infty_i$ for $i=1,2$.

\hspace*{-0.7 cm} \tabcolsep 0.04in
\begin{tabular}{llll}
$(8,9,0,1,\infty_1)$,& $(1,18,0,3,\infty_1)$,&
$(1,4,0,5,\infty_1)$,& $(1,6,0,7,\infty_1)$,\\
$(2,7,0,9,\infty_1)$,& $(2,11,0,1,\infty_2)$,&
$(6,13,0,3,\infty_2)$,& $(3,16,0,7,\infty_2)$,\\
$(3,8,0,15,\infty_2)$,& $(3,12,0,9,\infty_2)$,& $(0,15,4,9,14)$.
\end{tabular}
\\ \qed

\begin{Lemma}
\label{15:5} There exists an $HS(3,K_4^{(3)}+e;15,5)$.
\end{Lemma}
\proof  The required design is constructed on $\{0,1,\ldots,14\}$
with a hole $\{10,11,12,13,14\}$. All $89$ blocks are listed below.

\hspace*{-0.3 cm}\tabcolsep 0.04in
\begin{tabular}{llll}
$(0,9,1,8,6)$,& $(0,2,1,3,13)$,& $(0,5,1,4,10)$,& $(0,1,6,7,9)$,\\
$(10,11,0,1,14)$,& $(1,12,0,13,2)$,& $(2,4,0,6,9)$,&
$(0,2,5,7,10)$,\\ $(0,2,8,10,11)$,& $(0,2,9,11,10)$,&
$(0,2,12,14,9)$,& $(3,4,0,7,8)$,\\ $(0,6,3,5,9)$,& $(8,11,0,3,12)$,&
$(3,9,0,10,5)$,& $(0,3,13,14,8)$,\\ $(0,4,8,12,10)$,&
$(0,4,9,13,12)$,& $(0,10,4,14,5)$,& $(0,13,5,8,2)$,\\
$(0,9,5,12,6)$,& $(0,14,5,11,7)$,& $(0,14,6,8,5)$,&
$(0,10,6,12,7)$,\\
 $(0,13,6,11,5)$,& $(0,9,7,14,5)$,&
$(0,10,7,13,6)$,& $(0,12,7,11,6)$,\\
\end{tabular}

\hspace*{-0.3 cm}\tabcolsep 0.04in
\begin{tabular}{llll}
 $(2,4,1,7,9)$,&
$(2,6,1,5,11)$,& $(1,8,2,11,6)$,& $(1,2,9,10,6)$,\\
$(1,2,13,14,5)$,& $(1,4,3,6,10)$,& $(1,3,5,7,9)$,&
$(1,3,8,10,13)$,\\
$(1,3,9,11,5)$,& $(1,3,12,14,8)$,& $(1,8,4,13,6)$,&
$(1,9,4,12,7)$,\\
$(1,14,4,11,0)$,& $(1,12,5,8,11)$,& $(1,13,5,9,6)$,&
$(1,10,5,14,6)$,\\ $(1,14,6,9,8)$,& $(1,10,6,13,12)$,&
$(1,12,6,11,10)$,& $(1,8,7,14,12)$,\\ $(1,7,10,12,9)$,&
$(1,11,7,13,5)$,& $(2,5,3,4,8)$,& $(2,6,3,7,11)$,\\ $(2,3,8,9,14)$,&
$(10,11,2,3,14)$,& $(3,13,2,12,1)$,& $(8,14,2,4,9)$,\\
$(4,12,2,10,7)$,& $(2,4,11,13,9)$,& $(2,9,5,14,12)$,&
$(2,13,5,10,9)$,\\ $(2,5,11,12,9)$,& $(2,6,8,12,13)$,&
$(2,9,6,13,5)$,& $(2,6,10,14,8)$,\\ $(2,13,7,8,10)$,&
$(2,12,7,9,8)$,& $(2,7,11,14,8)$,& $(3,4,9,14,10)$,\\
$(3,4,10,13,9)$,& $(3,4,11,12,8)$,& $(3,14,5,8,7)$,&
$(3,12,5,10,8)$,\\ $(3,5,11,13,8)$,&
$(3,6,8,13,9)$,& $(3,6,9,12,8)$,& $(3,6,11,14,9)$,\\
 $(3,8,7,12,13)$,& $(3,13,7,9,11)$,&
$(3,14,7,10,11)$,& $(4,5,6,7,14)$,\\ $(4,5,8,9,11)$,&
$(4,11,5,10,6)$,& $(4,13,5,12,7)$,& $(4,8,6,10,7)$,\\
$(4,9,6,11,8)$,& $(4,12,6,14,13)$,& $(4,11,7,8,6)$,&
$(4,7,9,10,8)$,\\ $(4,7,13,14,9)$.
\end{tabular}
\\\qed

\begin{Lemma}
\label{16:6} There exists an $HS(3,K_4^{(3)}+e;16,6)$.
\end{Lemma}
\proof  The required design is constructed on $\{0,1,\ldots,15\}$
with a hole $\{10,11,12,13,14,15\}$. All $108$ blocks are listed
below.

\hspace*{-0.4 cm}\tabcolsep 0.04in
\begin{tabular}{llll}
$(0,12,7,11,10)$,& $(1,7,2,4,3)$,& $(1,2,5,6,11)$,&
$(1,2,8,11,5)$,\\
$(1,2,9,10,5)$,& $(1,2,12,15,8)$,& $(1,2,13,14,8)$,&
$(1,3,4,6,15)$,\\
$(1,7,3,5,4)$,& $(1,3,8,10,5)$,& $(1,3,9,11,5)$,& $(1,3,12,14,5)$,\\
$(1,3,13,15,7)$,& $(1,8,4,13,7)$,& $(1,9,4,12,14)$,&
$(1,15,4,10,7)$,\\ $(1,11,4,14,6)$,& $(1,8,5,12,7)$,&
$(1,13,5,9,7)$,& $(1,10,5,14,6)$,\\ $(1,11,5,15,7)$,&
$(1,15,6,8,7)$,& $(1,14,6,9,7)$,& $(1,10,6,13,7)$,\\
$(1,11,6,12,4)$,& $(1,14,7,8,5)$,& $(1,15,7,9,4)$,&
$(1,12,7,10,6)$,\\ $(1,11,7,13,5)$,& $(2,3,6,7,11)$,&
$(2,3,8,9,6)$,& $(2,3,10,11,6)$,\\
$(2,3,12,13,6)$,& $(2,3,14,15,8)$,& $(2,8,4,14,7)$,&
$(2,15,4,9,10)$,\\$(2,10,4,12,15)$,&
$(2,11,4,13,15)$,& $(2,15,5,8,6)$,& $(2,14,5,9,6)$,\\
$(2,10,5,13,6)$,& $(2,11,5,12,6)$,& $(2,8,6,12,7)$,&
$(2,9,6,13,14)$,\\ $(2,14,6,10,5)$,& $(2,15,6,11,8)$,&
$(2,8,7,13,12)$,& $(2,12,7,9,10)$,\\ $(2,10,7,15,12)$,&
$(2,11,7,14,5)$,& $(3,4,8,15,9)$,& $(3,4,9,14,8)$,\\
$(3,10,4,13,14)$,& $(0,1,14,15,7)$,& $(0,2,4,6,13)$,&
$(0,2,5,7,10)$,\\ $(0,2,8,10,9)$,&
$(0,2,9,11,8)$,& $(0,2,12,14,7)$,& $(0,2,13,15,6)$,\\
$(0,3,4,7,15)$,& $(0,3,5,6,15)$,& $(0,3,8,11,12)$,&
$(0,3,9,10,6)$,\\
$(0,3,12,15,6)$,& $(0,3,13,14,5)$,& $(0,4,8,12,14)$,&
$(0,4,9,13,14)$,\\
\end{tabular}

\hspace*{-0.4 cm}\tabcolsep 0.04in
\begin{tabular}{llll}
 $(0,4,10,14,9)$,& $(0,4,11,15,8)$,&
$(0,5,8,13,11)$,& $(0,5,9,12,15)$,\\
$(0,5,10,15,8)$,& $(0,5,11,14,9)$,& $(0,6,8,14,11)$,&
$(0,6,9,15,10)$,\\ $(0,6,10,12,8)$,& $(0,6,11,13,9)$,&
$(0,7,8,15,13)$,& $(0,7,9,14,12)$,\\ $(0,7,10,13,8)$,&
$(3,4,11,12,9)$,& $(3,5,8,14,10)$,& $(3,5,9,15,13)$,\\
$(3,10,5,12,15)$,& $(3,11,5,13,15)$,& $(3,6,8,13,9)$,&
$(3,6,9,12,10)$,\\ $(3,10,6,15,14)$,& $(3,11,6,14,12)$,&
$(3,8,7,12,4)$,& $(3,7,9,13,10)$,\\ $(3,10,7,14,13)$,&
$(3,15,7,11,5)$,& $(4,5,6,7,14)$,& $(4,5,8,9,12)$,\\
$(4,5,10,11,9)$,&
$(4,5,12,13,9)$,& $(4,5,14,15,9)$,& $(4,6,8,10,7)$,\\
$(4,6,9,11,15)$,& $(4,8,7,11,9)$,& $(0,1,2,3,5)$,& $(0,1,4,5,2)$,\\
$(0,1,6,7,15)$,& $(0,1,8,9,7)$,& $(0,1,10,11,8)$,& $(0,1,12,13,8)$.
\end{tabular}
\\\qed

\section{Conclusion}

In this section we give the necessary and sufficient conditions for
the existence of an $S(3,K_4^{(3)}+e,v)$.

\begin{Lemma}
\label{56}  There does not exist an $S(3,K_4^{(3)}+e,v)$ for
$v=5,6$.
\end{Lemma}
\proof Let $(X,{\cal B})$ be an $S(3,K_4^{(3)}+e,v)$ for $v=5,6$.
For any $x\in X$ denote by $d_i(x)$, $i=1,3,4$, the number of blocks
in $\cal B$ in which the degree of $x$ is $i$. It follows that
$d_1(x)+3d_3(x)+4d_4(x)=(v-1)(v-2)/2$. Note that
$d_1(x)+d_3(x)+d_4(x)\leq |{\cal B}|$.

For $v=5$, solving the equation with the constraint $|{\cal B}|=2$,
we have that for any $x\in X$, $d_1(x)=d_4(x)=0$ and $d_3(x)=2$. It
is easy to see that it is impossible. That is a contradiction.

For $v=6$, solving the equation with the constraint $|{\cal B}|=4$,
we have that:

($1$) $d_1(x)=0$, $d_3(x)=2$ and $d_4(x)=1$ (such vertex $x$ is
called a-vertex);

($2$) $d_1(x)=1$, $d_3(x)=3$ and $d_4(x)=0$ (such vertex $x$ is
called b-vertex);

($3$) $d_1(x)=2$, $d_3(x)=0$ and $d_4(x)=2$ (such vertex $x$ is
called c-vertex).

\noindent Denote by $\alpha$, $\beta$, $\gamma$ the number of
a-vertices, b-vertices and c-vertices respectively. Since each block
of ${\cal B}$ contains exactly one vertex with degree $1$, we have
that $\beta+2\gamma=|{\cal B}|=4$. Since each block of ${\cal B}$
contains exactly two vertices with degree $3$, we have that
$2\alpha+3\beta=2|{\cal B}|=8$. Due to $\alpha+\beta+\gamma=6$, we
have that $\alpha=4$, $\beta=0$ and $\gamma=2$. Let
$X=\{1,2,\ldots,6\}$, in which $1,2$ are c-vertices and $3,4,5,6$
are a-vertices. Without loss of generality we can assume that the
$4$ blocks of $\cal B$ are $(5,6,2,3,1)$, $(3,r,2,4,1)$,
$(3,*,1,5,2)$, $(*,*,1,6,2)$ or $(*,*,2,3,1)$, $(*,*,2,4,1)$,
$(3,4,1,5,2)$, $(3,r,1,6,2)$. Obviously any $r$ from
$\{1,2,\ldots,6\}$ is impossible. That is a contradiction. \qed

\begin{Lemma}
\label{10r+012}  There exists an $S(3,K_4^{(3)}+e,v)$ for $v\equiv
0,1,2$ $({\rm mod}$ $10)$.
\end{Lemma}
\proof By Lemmas \ref{71116} and \ref{101215}, there exists an
$S(3,K_4^{(3)}+e,v)$ for $v=10$, $11$, $12$, $31$, $32$. There
exists a $CS(3,K_4^{(3)}+e,30)$ of type $(15^2:0)$ from Lemma
\ref{g^n:0}. Then apply Construction \ref{filling construction} with
an $S(3,K_4^{(3)}+e,15)$ from Lemma \ref{101215} to obtain an
$S(3,K_4^{(3)}+e,30)$.

By Lemma \ref{GDD-46}, there exists a $GDD(3,\{4,6\},2n)$ of type
$2^n$ for any integer $4\leq n\leq 27$ and $n\neq 5,21$. Then apply
Construction \ref{GDD-FC} with a $GDD(3,K_4^{(3)}+e,5k)$ of type
$5^k$ for $k=4,6$ from Lemma \ref{GDD-5^4-5^6-10^5} to obtain a
$GDD(3,K_4^{(3)}+e,10n)$ of type $10^n$.

By Lemma \ref{45679} there exists an
$S(3,\{4,5,6,7,9,11,13,15,19,23$, $27\},u)$ for any integer $u\geq
4$, which implies the existence of a
$GDD(3,\{4,5,6,7,9,11,13,15,19,23,27\},u)$ of type $1^u$. Then apply
Construction \ref{GDD-FC} to obtain a $GDD(3,K_4^{(3)}+e,10u)$ of
type $10^u$, where the needed $GDD(3,K_4^{(3)}+e,50)$ of type $10^5$
is from Lemma \ref{GDD-5^4-5^6-10^5}, and the needed
$GDD(3,K_4^{(3)}+e,10m)$ of type $10^m$ for
$m=4,6,7,9,11,13,15,19,23,27$ have been constructed in the second
paragraph.

Start with the resulting $GDD(3,K_4^{(3)}+e,10u)$ of type $10^u$,
and apply Construction \ref{g^2-construction} to obtain a
$CS(3,K_4^{(3)}+e,10u+s)$ of type $(10^u:s)$ for $u\geq 2$, $u\neq
3$ and $s=0,1,2$, where the needed $CS(3,K_4^{(3)}+e,20+s)$ of type
$(10^2:s)$ are from Lemmas \ref{g^n:0}, \ref{g^n:1} and \ref{g^n:2}.
Fill in holes by Construction \ref{filling construction} to obtain
an $S(3,K_4^{(3)}+e,10u+s)$. \qed

\begin{Lemma}
\label{10r+56}  There exists an $S(3,K_4^{(3)}+e,v)$ for $v\equiv
5,6$ $({\rm mod}$ $10)$ and $v\geq 15$.
\end{Lemma}
\proof By Lemmas \ref{71116} and \ref{101215}, there exists an
$S(3,K_4^{(3)}+e,v)$ for $v=15$, $16$, $26$. There exists a
$CS(3,K_4^{(3)}+e,25)$ of type $(6^4:1)$ from Lemma \ref{g^n:1}.
Then apply Construction \ref{filling construction} with an
$S(3,K_4^{(3)}+e,7)$ from Lemma \ref{71116} to obtain an
$S(3,K_4^{(3)}+e,25)$.

By Lemma \ref{2^n:2}, there exists a $CS(3,\{4,6\},2u+2)$ of type
$(2^u:2)$ for any integer $u\geq 3$. Then applying Construction
\ref{HFC} with $b=5$ and $r=0,1$, we obtain a
$CS(3,K_4^{(3)}+e,10u+5+r)$ of type $(10^u:5+r)$, where the needed
$CS(3,K_4^{(3)}+e,5k+r-5)$ of type $(5^{k-1}:r)$ for $k=4,6$ are
from Lemmas \ref{g^n:0} and \ref{g^n:1}, and the needed
$GDD(3,K_4^{(3)}+e,5k)$ of type $5^k$ for $k=4,6$ are from Lemma
\ref{GDD-5^4-5^6-10^5}. By Lemmas \ref{15:5} and \ref{16:6}, there
exist an $HS(3,K_4^{(3)}+e;15,5)$ and an $HS(3,K_4^{(3)}+e;16,6)$.
Thus apply Construction \ref{filling construction} to obtain an
$S(3,K_4^{(3)}+e,10u+5+r)$. \qed

\begin{Lemma}
\label{10r+7}  There exists an $S(3,K_4^{(3)}+e,v)$ for $v\equiv 7$
$({\rm mod}$ $10)$.
\end{Lemma}

\proof By Lemma \ref{46}, there exists an $S(3,\{4,6\},u+1)$
$(X,{\cal B})$ with $u\equiv 1$ $({\rm mod}$ $2)$ and $u\geq 3$.
Then for any $x_0\in X$, $(X,\{x_0\},\{\{x\} :x\in X\setminus
\{x_0\}\},{\cal B})$ is a $CS(3,\{4,6\},u+1)$ of type $(1^u:1)$.
Applying Construction \ref{HFC} with $b=5$ and $r=2$, we have a
$CS(3,K_4^{(3)}+e,5u+2)$ of type $(5^u:2)$, where the needed
$CS(3,K_4^{(3)}+e,5k-3)$ of type $(5^{k-1}:2)$ for $k=4,6$ are from
Lemma \ref{g^n:2}, and the needed $GDD(3,K_4^{(3)}+e,5k)$ of type
$5^k$ for $k=4,6$ are from Lemma \ref{GDD-5^4-5^6-10^5}. By Lemmas
\ref{71116} there exists an $S(3,K_4^{(3)}+e,7)$. Thus apply
Construction \ref{filling construction} to obtain an
$S(3,K_4^{(3)}+e,5u+2)$. \qed

\noindent {\bf Proof of Theorem \ref{K_4+e}}: The necessity follows
from Lemmas \ref{necessary} and \ref{56}. The sufficiency follows
from Lemmas \ref{10r+012}-\ref{10r+7}. \qed


\begin{thebibliography}{Z}

\bibitem{berge}
C. Berge, {\it Hypergraphs: Combinatorics of Finite Sets}, North
Holland, Amsterdam, 1989.

\bibitem{bjl}
T. Beth, D. Jungnickel and H. Lenz, {\it Design Theory}, Cambridge
University Press, Cambridge, UK, 1999.

\bibitem{be} D. Bryant, S. El-Zanati, Graph decompositions,
in: {\it CRC Handbook of Combinatorial Designs} (C. J. Colbourn and
J. H. Dinitz eds.) CRC Press, 2006, 477-487.

\bibitem{cgal}
G. D. Crescenzo and C. Galdi, {\it Hypergraph decomposition and
secrete sharing}, Proceedings of ISAAC'$03$, LNCS $2906$,
Springer-Verlag, (2003), 645-654.

\bibitem{fc}
T. Feng and Y. Chang {\it Decompositions of the $3$-uniform
hypergraphs $K_v^{(3)}$ into hypergraphs of a certain type}, Science
in China(A), 50(2007), 1035-1044.

\bibitem{hanani60}
H. Hanani, {\it On quadruple systems}, Canad. J. Math., 12(1960),
145-157.

\bibitem{hanani63}
H. Hanani, {\it On some tactical configurations}, Canad. J. Math.,
15(1963), 702-722.

\bibitem{ji45}
L. Ji, {\it On the $3$BD-closed set $B_3(\{4,5\})$}, Discrete Math.,
287(2004), 55-67.

\bibitem{ji}
L. Ji, {\it More results on perfect $(3,3,6k+4;6k-2)$-threshold
schemes}, J. Combin. Des., 15(2007), 151-166.

\bibitem{m74}
W. H. Mills, {\it On the covering of triples by quadruples}, Congr.
Numer., 10(1974), 563-586.

\bibitem{mr}
H. Moh$\acute{a}$csy and D. K. Ray-Chaudhuri, {\it Candelabra
systems and designs}, J. Statist. Plann. Inference., 106(2002),
419-448.

\end{thebibliography}
\end{document}